\documentclass[twoside]{irmaems}
\usepackage{amssymb} 
\usepackage{amsmath} 
\usepackage{latexsym}

\renewcommand\theequation{\arabic{section}.\arabic{equation}}

\theoremstyle{definition} 

\theoremstyle{plain}      

\input epsf

\begin{document}

\title{Variational Principles on Triangulated Surfaces}

\author{Feng Luo\thanks{ Work partially supported by the NSF.} }

\address{
Department of Mathematics\\
Rutgers University\\
Piscataway, NJ 08854, USA\\
email:\,\tt{ fluo\@math.rutgers.edu}\\
\\
\\
\text{To Professor S.-T. Yau on his sixtieth birthday} }
\maketitle

\begin{abstract} We give a brief introduction to some of the
recent works on finding geometric structures on triangulated
surfaces using variational principles.
\end{abstract}
\tableofcontents

\section{Introduction}\label{s-1}

The main objects of  study are geometry of polyhedral surfaces.
Take a finite set of points in the Euclidean 3-space
$\mathbb{E}^3$ and the convex hull of these points. We obtain a
convex polytope whose vertices are among the given finite set. If
the points are generic, then the convex polytope has triangle
faces. In this case, the boundary surface is a polyhedral surface.
It has two properties. First the surface is triangulated and
second the induced metric on the surface is locally flat except at
the vertices. Recall that a triangulation of a surface is defined
as follows. Take a finite collection of disjoint triangles and
identify pairs of edges by homeomorphisms. The quotient space is a
surface with a triangulation whose cells are the quotients of
triangles, edges and vertices in the disjoint union.
\\
\\
{\bf Definition 1.1.} A \it Euclidean polyhedral surface \rm is a
triple $(S,T,d)$ where $S$ is a closed surface, $T$ is a
triangulation of $S$, and $(S, d)$ is a metric space with metric
$d$ so that the restriction of $d$ to each triangle is isometric
to a Euclidean triangle. We will call the metric $d$ a \it
Euclidean polyhedral metric. \rm The \it discrete curvature \rm
(or simply curvature) $k_0$ of $(S,T,d)$ is a function which
assigns each vertex $2\pi$ less the sum of inner angles at the
vertex, i.e.,
$$k_0(v)=2\pi - \sum_{i=1}^m \theta_i,$$
where $\theta_1,...,\theta_m$ are the inner angles (of triangles)
at the vertex $v$. See Figure 3.1.
\medskip

We will also call above a $\mathbb{E}^2$ polyhedral surface. If we
use the spherical (or hyperbolic) triangles instead of Euclidean
triangles in definition 1.1, the result is called a spherical (or
$\mathbb{S}^2$) polyhedral surface (resp. $\mathbb{H}^2$ or
hyperbolic polyhedral surface). Spherical and hyperbolic
polyhedral surfaces have been studied extensively.  The discrete
curvature $k_0$ is defined by the same formula for $\mathbb{S}^2$
and $\mathbb{H}^2$ polyhedral surfaces.

From the definition, it is clear that the basic unit of discrete
curvature is the inner angle. Furthermore, the metric-curvature
relation is given by the cosine law. Just like the smooth case,
one of the main problems of study in polyhedral surface is to
understand the relationship between the metric and its curvature.
Naturally we should study the cosine law carefully.

The goal of the paper is to introduce some of the recent works on
finding geometric structures on triangulated surfaces using
variational principles.
 Ever since W. Thuston's work on geometrization of
Haken 3-manifolds and circle packing in 1978, there have been many
works in this area.
 The key step in
the variational framework is to define an appropriate action
functional so that the critical points of the functional are the
geometric structures that one is seeking.

The first such action functional was discovered in a seminal work
by Colin de Verdi\`ere \cite{CV} in 1991 for circle packing
metrics. Colin de Verdi\`ere's functional is derived from the
Schlaefli formula for volume of tetrahedra. In \cite{CV}, Colin de
Verdi\`ere introduced the first variational principle on
triangulated surfaces in recent times and gave a proof of
Thurston-Andreev's existence and uniqueness theorem for circle
packing metrics. In a remarkable paper \cite{Ri} in 1994, I. Rivin
used the 3-dimensional volume of a hyperbolic ideal tetrahedra as
the action functional and established a beautiful variational
principle for Euclidean polyhedral surfaces. Since then, many
other variational principles for polyhedral surfaces have been
established. See for instance, \cite{Br}, \cite{BS}, \cite{CKP},
\cite{Le}, \cite{Lu1}, \cite{Lu2}, \cite{Sch}, \cite{Sp} and
others. Amazingly, almost all action functionals discovered so far
are related to the Schlaefli formula. The only exception is in the
beautiful work of \cite{BS}. The action functional is derived from
a discrete integrable system. Very recently, we realized
\cite{Lu3} that the cosine law and its derivative form are rich
sources for action functionals and these include all the previous
approaches. These recently discovered action functionals, when
view from some perspectives, can be considered as 2-dimensional
counterparts of the Schlaefli formula. The complete list of all
possible 2-dimensional counterparts of the Schlaefli formula has
been found in \cite{Lu3}.

This paper is organized as follows. In section 2, we discuss the
construction of action functionals in 2-dimension. In section 3,
we discuss various variational principles associated to the action
functionals.  In section 4, we discuss the problem of moduli space
of all curvatures. The last section addresses some open problems.

\section{The Schlaefli formula and its counterparts in dimension 2}

One of the most beautiful identities in low-dimensional geometry
is the Schlaefli formula. It states that for a tetrahedron in a
constant curvature $\lambda =\pm 1$ space, the volume $V$, the
length $l_i$ and the dihedral angle $a_i$ at the $i$-th edge are
related by

\begin{equation}
\label{2.1}  \frac{\partial V}{\partial a_i} = \frac{\lambda}{2}
l_i
\end{equation}
 where $V=V(a_1, ..., a_6)$ is  a
function of the angles. The formula relates the three most
important 3-dimensional geometric quantities: volume, metric
(=length), and curvature (=dihedral angles) in a simple identity.
One consequence of (2.1) is that the differential 1-forms,
\begin{equation}\label{2.2}
\sum_{i=1}^3 l_i da_i \quad  \text{and} \quad \sum_{i=1}^3 a_i
dl_i  \quad \text{ are closed.} \end{equation} Indeed, $2dV =
\lambda \sum_{i=1}^3 l_i da_i$. One can recover (2.1) from (2.2)
by integration. By taking the Legendre transformation, one obtains
$ H(l_1, ..., l_6) = \sum_{i=1}^6 l_i a_i - 2 \lambda V$ so that
\begin{equation} \label{2.3} \frac{\partial H}{\partial l_i} =
a_i. \end{equation}

Identities (2.1) and (2.3) are the starting points of several
variational principles for finding constant curvature metrics on
triangulated 3-manifolds. They are the basic ingredients in Regge
calculus in physics which is a discretized general relativity.

\subsection{Regge calculus and Casson's approach in dimension 3}
Here is an illustration of the use of (2.3) after A. Casson and
others. Fix a triangulated closed 3-manifold $(M, T)$ and consider
the space $X$ of all hyperbolic polyhedral metrics on $(M, T)$. By
definition, a hyperbolic polyhedral metric on $(M, T)$ is a metric
on $M$ so that the restriction of the metric to each tetrahedron
is isometric to a hyperbolic tetrahedron. Let $E$ be the set of
edges in the triangulation $T$. Then a polyhedral metric on $(M,
T)$ is determined by the edge length function $l : E \to
\mathbf{R}$ sending an edge to its length. The \it discrete
curvature \rm (or simply curvature)  $K$ of the metric $l$ is the
function $K: E \to \mathbf{R}$ sending an edge to $2\pi$ less the
sum of dihedral angles at the edge. If the curvature $K=0$, then
the polyhedral metric is a smooth hyperbolic metric. This is
proved as follows. First, the curvature $K=0$ shows that the
metric is smooth at the interior points of each edge. To show that
it is also smooth at each vertex, one considers the spherical link
at the vertex. The link is isometric to a spherical polyhedral
2-sphere with discrete curvature $k_0=0$ at each vertex. This
shows the link is isometric to the standard 2-sphere. Thus the
metric on the 3-manifold is smooth at each vertex. Now in Casson's
approach, one defines the Einstein-Hilbert action of a polyhedral
metric $l$ to be:
$$F(l)= -\sum_{ \sigma_i} H(l_{i_1}, ..., l_{i_6}) +  2\pi \sum_{e_j} l_j$$ where the first
sum is over all tetrahedra $\sigma_i$ with six edge lengths
$l_{i_1}, ..., l_{i_6}$ and the second sum is all edges $e_j$ of
length $l_j$. One easy consequence of (2.3) is that the
Euler-Lagrangian equation for $F$, considered as a function
defined on $X$, is given by
\begin{equation}\label{2.4}
 \frac{\partial F}{\partial l_i} = K_i  \end{equation} where $l_i$
and $K_i$ are the length and curvature at the $i$-th edge. In
particular, it follows that the critical points of the
Einstein-Hilbert functional are the hyperbolic metrics.  To derive
the Euler-Lagrangian equation (2.4), let us assume that the
dihedral angles at the $i$-th edge are $a_1, ..., a_n$ and the
tetrahedra adjacent to the $i$-th edge are $\sigma_1, ...,
\sigma_n$.  Then by the definition of $F$ and the Schlaefli
formula (2.3) applied to each $\sigma_j$, we have
$$ \frac{\partial F}{\partial l_i} =(-a_1 - ...-a_n) + 2\pi =K_i.$$

The above approach is ubiquitous in variational framework on
triangulated spaces.  The important feature of this approach is
that the action functional is local. This means that the value of
the action functional on a polyhedral metric is the sum of the
functional on its top dimensional simplexes. Thus the main issue
for variational framework for triangulated surfaces is to find
action functionals for geometric triangles.

\subsection{The work of Colin de Verdi\`ere, Rivin,
Cohen-Kenyon-Propp and Leibon} For a long time, the Gauss-Bonnet
formula for area of triangles had been considered as the only
counterpart of the Schlaefli formula in dimension 2. This is
probably due to the view point that one should emphasize the role
of volume in (2.1). The view changed when Colin de Verdi\`ere
\cite{CV} produced the first striking 2-dimensional counterpart of
the (1.1) by paying attention to (2.2). In this section, we will
introduce briefly the action functionals in the work of \cite{CV},
\cite{Ri}, \cite{CKP}, and \cite{Le}.

In the work of \cite{CV} which gives a new proof of
Andreev-Thurston's existence and uniqueness of circle packing
metrics using variational principle, Colin de Verdi\`ere considers
triangles of edge lengths $r_1+r_2, r_2+r_3, r_3+r_1$ and angles
$a_1, a_2, a_3$ where the angle $a_i$ faces the edge of length
$r_j+r_k$, \{i,j,k\}=\{1,2,3\}. See Figure 2.1 (a). He proved that
the following three 1-forms $w$ are closed. These are the
counterparts of (2.2) in 2-dimension.

For a Euclidean triangle, the 1-form \begin{equation} \label{2.5}
w =\sum_{i=1}^3 \frac{a_i}{r_i} dr_i = \sum_{i=1}^3 a_i d \ln r_i
\end{equation} is closed. For a hyperbolic triangle, the 1-form
\begin{equation} \label{2.6} w =\sum_{i=1}^3 \frac{a_i}{\sinh r_i}
d r_i = \sum_{i=1}^3 a_i d \ln \tanh(r_i/2) \end{equation} is
closed. For a spherical triangle, the 1-form \begin{equation}
\label{2.7} w =\sum_{i=1}^3 \frac{ a_i}{\sin r_i} d r_i =
\sum_{i=1}^3 a_i d \ln \tan(r_i/2) \end{equation} is closed.
Furthermore, write the 1-form $w $ as $\sum_{i=1}^3 a_i du_i$.
Then he proved that the integration $F(u) = \int^u w$ is concave
in $u =(u_1, u_2, u_3)$ in the cases of Euclidean and hyperbolic
triangles. We are informed by Colin de Verdi\`ere that these
1-forms were discovered by considering the Schlaefli formula.

\medskip

\epsfxsize=3.0truein\centerline{\epsfbox{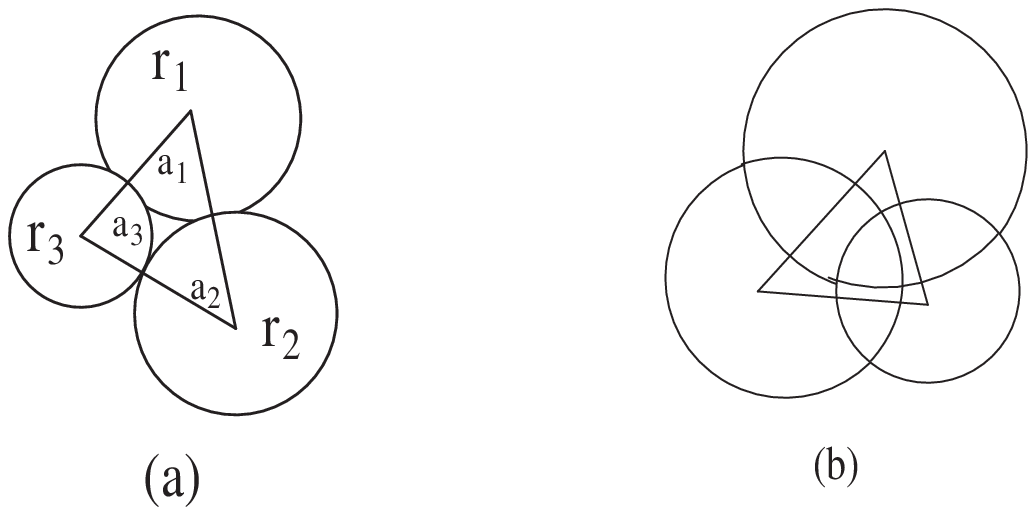}}

\centerline{Figure 2.1}



The Legendre transformation of $F(u)$ for Euclidean triangle was
discussed in the work of Braegger \cite{Br} who identifies it as
the volume of a hyperbolic tetrahedron. The geometric meaning of
other $F(u)$'s or its Legendre transformation is not known.

Thurston's original work on circle packing allows circles to
intersect at angles. See Figure 2.1(b). For these intersecting
angle circle packing metrics, Colin de Verdi\`ere's 1-form are
still closed. This is proved in \cite{CL}.

In the work of Rivin \cite{Ri} and Cohen-Kenyon-Propp \cite{CKP},
they consider Euclidean triangles of edge lengths $l_1, l_2, l_3$
so that the opposite angles are $a_1, a_2, a_3$. The action
functional considered by Rivin is the Lengendre transformation of
that of Cohen-Kenyon-Propp.  It is proved that the  1-form
\begin{equation} \label{2.8} w =\sum_{i=1}^3 \frac{a_i}{l_i} dl_i = \sum_{i=1}^3 a_i d
\ln l_i  \end{equation} is closed. The integral $F(u) = \int^u w$
is shown to be convex in $u=(\ln l_1, \ln l_2, \ln l_3)$ (see
\cite{Ri}). It is proved in \cite{Ri} that the Legendre
transformation of $F(u)$ is the volume of the hyperbolic ideal
tetrahedron with dihedral angles $a_1, a_2, a_3, a_1, a_2, a_3$.

In the work of \cite{Le}, Leibon considers a hyperbolic triangle
of edge lengths $l_i$ so that the opposite angle is $r_j+r_k$
where $\{i,j,k\}=\{1,2,3\}$. He proved that the 1-form
\begin{equation} \label{2.9} w= \sum_{i=1}^3 \ln \sinh(l_i/2) dr_i \end{equation}
 is closed and its
integration is strictly convex in $(r_1, r_2, r_3)$.  Furthermore,
the integration $\int^r w$ is proved in \cite{Le} to be the volume
of a hyperbolic prism.

Evidently, (2.5)-(2.9) should be considered as 2-dimensional
counterparts of the Schlaefli formula (2.2). We will discuss the
applications of these action functionals in section 3.

\subsection{The Cosine Law and 2-dimensional Schlaefli formulas}
The 2-dimensional Schlaefli formulas that we are seeking are some
relationship between the lengths and angles of a triangle. Let
$\mathbf{E}^2, \mathbf{H^2}$ and $\mathbf{S^2}$ be the Euclidean
plane, the hyperbolic plane and the 2-sphere respectively. Given a
triangle in $\bold H^2$, $\bold E^2$ or $\bold S^2$ of inner
angles $\theta_1, \theta_2, \theta_3$ and edge lengths $l_1, l_2,
l_3$ so that $\theta_i$ is facing the $l_i$-th edge, the cosine
law expressing length $l_i$ in terms of the angles $\theta_r$'s
is,
$$ \cos(\sqrt{ \lambda  }l_i) =
 \frac{ \cos \theta_i + \cos \theta_j \cos \theta_k}{\sin \theta_j \sin \theta_k}$$
where $ \lambda  =1,-1,0$ is the curvature of the space $\bold
S^2$, or $\bold H^2$ or $\bold E^2$ and $\{i,j,k\}=\{1,2,3\}$.
 Another related cosine law
is$$ \cosh(l_i) = \frac{ \cosh \theta_i + \cosh \theta_j \cosh
\theta_k}{\sinh \theta_j \sinh \theta_k}$$ for a right-angled
hyperbolic hexagon with three non-adjacent edge lengths $l_1, l_2,
l_3$ and their opposite edge lengths $\theta_1, \theta_2,
\theta_3$.

These formulas suggest that we should consider the following.
Suppose a function $y=y(x)$ where $y=(y_1, y_2, y_3) \in \bold
C^3$  and $x=(x_1, x_2, x_3)$ is in some open connected set in
$\bold C^3$ so that $x_i$'s and $y_i$'s are related by

\begin{equation}\label{2.10}
 \cos ( y_i) = \frac{\cos x_i + \cos x_j \cos x_k}{\sin (x_j)
\sin (x_k)} \end{equation} where $\{ i,j,k \}=\{ 1,2,3 \}$. We say
$y=y(x)$ is the \it  cosine law function\rm. Let $r_i
=\frac{1}{2}(x_j+x_k-x_i) $. Then $r=(r_1, r_2, r_3)$ is a new
parametrization so that $x_i = r_j + r_k$. We will also consider
the composition function $y=y(r_1, r_2, r_3)$.

The following proposition establishes the basic properties of the
cosine law function. The proof of the proposition is a simple
exercise in calculus and will be omitted. See \cite{Lu3} for
details.

\medskip
\noindent {\bf Proposition 2.1.} (\cite{Lu3}) \it Suppose the
cosine law function $y=y(x)$ is defined on an open connected set
in $\bold C^3$ which contains a point $(a,a,a)$ so that $y(a,a,a)
= (b,b,b)$. Let $ A_{ijk} = \sin y_i \sin x_j \sin x_k$ where
$\{i,j,k\}=\{1,2,3\}$. Then \begin{equation}\label{2.11} A_{ijk} =
A_{jki}, \quad i.e., \quad
\frac{\sin(y_i)}{\sin(x_i)}=\frac{\sin(y_j)}{\sin(x_j)}
\end{equation}

\begin{equation} \label{2.12}
\frac{\tan(y_i/2)}{\cos(r_i)}=\frac{\tan(y_j/2)}{\cos(r_j)}
\end{equation}

$$A^2_{ijk} =  1- \cos^2 x_i -\cos^2 x_j -\cos^2 x_k - 2 \cos x_i \cos x_j \cos x_k. $$
At a point $x$ where $A_{ijk} \neq 0$, then,
$$ \frac{ \partial y_i}{\partial x_i }= \frac{\sin x_i}{A_{ijk}},    $$
$$ \frac{  \partial y_i}{\partial x_j} = \frac{\partial y_i}{\partial x_i} \cos y_k, $$
\begin{equation}\label{2.13}  \cos(x_i) = \frac{ \cos y_i - \cos y_j \cos y_k}{\sin y_j
\sin y_k}. \end{equation} \rm

In particular, (2.13) shows that the roles of $x$ and $y$ in the
cosine law are essentially symmetric. The identities (2.11) and
(2.12) are called the \it Sine Law \rm and the \it Tangent Law \rm
of the cosine law function.

The problem of finding all 2-dimensional counterparts of the
Schlaefli formula can be rephrased as follows. It corresponds to
generalizing Schlaefli identity (2.2).

\medskip
\noindent {\bf Problem 2.2}.  \it Suppose the cosine law function
$y=y(x)$ is defined on an open connected set in $\bold C^3$.

\noindent (i) Find all smooth non-constant functions $f(t)$ and
$g(t)$ so that the differential 1-form  $ w =\sum_{i=1}^3
f(y_i)dg(x_i)$  is closed.

\noindent (ii) Find all smooth non-constant functions $f(t)$ and
$g(t)$ so that the differential  1-form $ \sum_{i=1}^3 f(y_i)
dg(r_i)$ is closed where  $r_i = \frac{x_j+x_k-x_i}{2}$,
$\{i,j,k\}=\{1,2,3\}$.

\rm
\medskip

If we find these 1-forms $w$, then the integrals $F = \int w$ will
be used as action functionals for variational principles on
surfaces. By the construction of the 1-form, the function $F$
satisfies
$$ \frac{ \partial F}{\partial g(x_i)} = f(y_i)$$
or $$\frac{\partial F}{\partial g(r_i)} = f(y_i).$$ These are all
2-dimensional counterparts of the Schlaefli formula (2.2).

Problem 2.2 was solved in \cite{Lu3}.

\medskip
\noindent {\bf Theorem 2.3. (\cite{Lu3})} \it The following is the
complete list of functions $f$ and $g$ up to scaling and complex
conjugation. There exists a complex number $h$ so that, in the
case (i), $$ f(t) = \int^t \sin^h(s) ds \quad \text{and} \quad
g(t) = \int^t \sin^{-h-1}(s) ds$$ and in the case (ii),
$$ f(t) = \int^t  \tan^h(s/2) ds \quad \text{and} \quad g(t) =
\int^t \cos^{-h-1}(s) ds.$$
 In particular, all
closed 1-forms $\sum_{i=1}^3 f(y_i) d g(x_i)$ and $\sum_{i=1}^3
f(y_i) d g(r_i)$ are holomorphic or anti-holomorphic. \rm

\medskip
The details of the proof the theorem can be found in \cite{Lu3}.
We give a sketch of the proof here. To verify that the 1-forms
listed above are closed is a straight forward calculation using
the sine law  $\frac{\sin (y_i)}{\sin (x_i)} =\frac{\sin
(y_j)}{\sin (x_j)}$ and the tangent law  $\frac{ \tan(y_i/2)}{\cos
(r_i)} =\frac{\tan(y_j/2)}{\cos(r_j)}$. The proof that these are
the complete list of all functions $f,g$ up to scaling and complex
conjugation is due to the uniqueness of the sine law and tangent
law. To be more precise, in the case of the sine law, it can be
shown that (\cite{Lu3}, lemma 2.3) if $f, g$ are two smooth
non-constant functions so that $\frac{f(y_i)}{g(x_i)} = \frac{
f(y_j)}{g(x_j)}$ for all $x$, then there are constants $ \lambda
$, $\mu$, $c_1, c_2$ so that, $f(t) =c_1\sin^{\lambda }(t)
\sin^{\mu}(\bar t)$ and $g(t) =c_2\sin^{\lambda}(t)\sin^{\mu}(\bar
t)$.

By specializing theorem 2.3 to triangles in $\bold S^2$, $\bold
E^2$ and $\bold H^2$ and integrating the 1-forms, we obtain
various energy functionals for variational framework on
triangulated surfaces.
 We have identified all those convex
or concave energies constructed in this way.

\medskip

\noindent {\bf Theorem 2.4 (\cite{Lu3}).}  \it  Let $l=(l_1, l_2,
l_3)$ and $\theta=(\theta_1, \theta_2, \theta_3)$ be lengths and
angles of a triangle in $\bold E^2$, $\bold H^2$ or $\bold S^2$.
Let $ h \in  \bold R$ and $u=(u_1, u_2, u_3)$.

The following is the complete list, up to scaling, of all closed
real-valued 1-forms of the form $\sum_{i=1}^3 f(\theta_i) d
g(l_i)$  for some non-constant smooth functions $f, g$ so that its
integral is convex or concave.

\noindent (i) For a Euclidean triangle,

$$ w_{ h  } = \sum_{i=1}^3 \frac{\int^{\theta_i} \sin^{ h  }(t) dt}{l_{i}^{ h  +1}}  dl_i.$$
Its  integral  $\int^u w_{ h  }$ is locally convex in variable $u$
where $u_i =\int^{l_i}_1 t^{- h   -1} dt$.

\noindent (ii) For a spherical triangle,
$$ w_{ h  } = \sum_{i=1}^3 \frac{\int^{\theta_i} \sin^{ h  }(t) dt}{\sin^{ h  +1}(l_i)} dl_i.$$
The integral $\int^u w_{ h  }$ is locally strictly convex in $u$
where
 $u_i =\int^{l_i}_{\pi/2}$  $ \sin^{- h  -1}(t) dt$.

The following are the complete list, up to scaling, of all closed
real-valued 1-forms of the form $\sum_{i=1}^3 f(l_i) d g(r_i)$
(where $\theta_i=r_j+r_k$) or $\sum_{i=1}^3 f(\theta_i) d g(r_i)$
(where $l_i=r_j+r_k$) for some non-constant smooth functions $f,g$
so that its integral is either convex or concave.

\noindent (iii) For a Euclidean triangle of angles $\theta_i$ and
opposite edge lengths $r_j+r_k$,
$$\eta_{ h  } =\sum_{i=1}^3 \frac{\int^{\theta_i} \cot^{ h  }(t/2) dt}{ r_i^{ h  +1}}
dr_i.$$  Its  integral  $\int^{u} \eta_{ h  }$ is locally concave
in $u=(u_1, u_2, u_3)$ where $u_i =\int_1^{r_i} t^{- h  -1} dt$.

\noindent (iv) For a hyperbolic triangle of angles $\theta_i$ and
opposite edge lengths $r_j+r_k$,
$$\eta_{ h  } = \sum_{i=1}^3 \frac{\int^{\theta_i} \cot^{ h   }(t/2) dt }{\sinh^{ h  +1}(r_i)} dr_i.$$
 Its  integral  $\int^{u}
\eta_{ h  }$ is locally strictly concave in $u$ where $u_i
=\int_1^{r_i} \sinh^{ -h-1 }(t) dt$.

\noindent (v) For a hyperbolic triangle of edge lengths $l_i$ and
opposite angles $r_j+r_k$,
$$\eta_{ h  } = \sum_{i=1}^3 \frac{\int^{l_i} \tanh^{ h  }(t/2) dt}{ \cos^{ h  +1}(r_i)} dr_i.$$
 Its  integral  $\int^{u}
\eta_{ h  }$ is locally strictly convex in $u$ where $u_i
=\int_1^{r_i} \cos^{ -h  -1}(t) dt$.

\noindent (vi) For
 a hyperbolic right-angled hexagon of three non-pairwise adjacent edge lengths $l_i$ and opposite edge lengths $r_j+r_k$,
$$\eta_{ h  } = \sum_{i=1}^3 \frac{\int^{l_i} \coth^{ h  }(t/2) dt }{\cosh^{ h  +1}(r_i)} dr_i.$$
 Its  integral  $\int^{u}
\eta_{ h  }$ is locally strictly concave in $u$ where $u_i
=\int_1^{r_i} \cosh^{ -h  -1}(t) dt$.\rm

\rm

\medskip

A sketch of the proof of theorem 2.4 goes as follows. First, by
theorem 2.3, these 1-forms are closed. It remains to show that the
integrals are convex or concave. This follows by showing that the
Hessian matrix of the integral is positive (or negative) definite.
We first observe that the Hessian matrices of $\int w_h$ (or $\int
\eta_h$) and $\int w_{h'}$ (or $\int \eta_{h'}$) are congruent for
different $h, h'$. Thus, it suffices to check the convexity for
one $h$. This was been achieved in various cases by different
authors. Indeed, case (i) for $h=0$ is proved in \cite{Ri}, case
(ii) for $h=0$ is proved in \cite{Lu1}, (iii) and (iv) for $h=0$
are proved in \cite{CV}, case (v) for $h=-1$ is proved in
\cite{Le}, and case (vi) for $h=-1$ is proved in \cite{Lu2}.

\subsection{The geometric meaning of some action functionals}
The geometric meaning of the integrals $\int w_h$ or $\int \eta_h$
are not known except in the following cases.  The Legendre
transform of $\int w_0$ for Euclidean triangles is the hyperbolic
volume of an ideal tetrahedron first discovered by Rivin. Leibon
showed that the integral $\int \eta_{-1}$ for hyperbolic triangle
is the volume of an ideal hyperbolic prism. We showed in
\cite{Lu1} that $\int w_0$ for a spherical triangle is the volume
of an ideal hyperbolic octahedron.  In these three cases, there is
a common way to describe the action functional. Given a Euclidean,
or a hyperbolic or a spherical triangle, we consider the triangle
to be drawn in the sphere at infinity of the hyperbolic 3-space
bounded by three circles. These three circles will intersect at a
finite set $X$ of points in $\mathbf{S^2}$ where $|X|=4$ for
Euclidean triangle and $|X| =6$ otherwise. Then the action
functional associated to the triangle is the hyperbolic volume of
the convex hull of $X$ in the hyperbolic 3-space. See Figure 2.2.

\medskip

\epsfxsize=3.5truein\centerline{\epsfbox{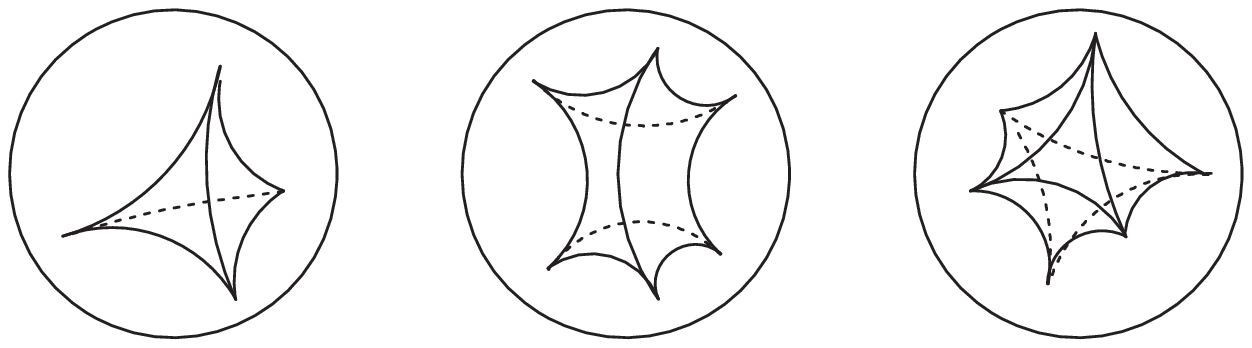}}

\centerline{Figure 2.2}
\medskip
\section{Variational principles on surfaces}
We discuss some  applications of the 2-dimensional Schlaefli
formulas in this section. The most prominent applications are
rigidity theorems for polyhedral surfaces. For simplicity, let us
assume that surfaces are closed in this section. All results can
be generalized without difficulties to compact surfaces with
boundary.

All of these applications to rigidity are based on the following
simple lemma.

\medskip
\noindent {\bf Lemma 3.1.} \it Suppose $\Omega \subset \bold R^n$
is an open convex set and $W: \Omega \to \bold R^n$ is a smooth
function with positive definite Hessian matrices. Then the
gradient $\bigtriangledown W: \Omega \to \bold R^n$ is a smooth
embedding. If $\Omega$ is only assumed to be open in $\bold R^n$,
then $\bigtriangledown W: \Omega \to \bold R^n$ is a local
diffeomorphism. \rm

\medskip
Indeed, consider the graph $G$ of the function $W$ in
$\mathbf{R^{n+1}}$. The convexity of $W$ shows that the graph $G$
is a strictly convex hypersurface. Thus normal vectors to the
graph $G$ at different points are not parallel. However, normal
vectors are of the form $(\bigtriangledown W, -1)$. It follows
that $\bigtriangledown W: \Omega \to \bold R^n$ is injective. The
Jacob matrix of $\bigtriangledown W$ is the Hessian of $W$. Thus
the map $\bigtriangledown W$ is a smooth embedding.

\medskip
\subsection{Colin de Verdi\`ere's proof of Thurston-Andreev's
rigidity theorem}

In his work \cite{Th} on constructing hyperbolic metrics on
3-manifolds, Thurston introduced circle packing  metrics on a
triangulated surface $(S, T)$.  Let $V$ and $E$ be the sets of
vertices and edges in the triangulation $T$.
 A \it circle packing metric \rm on $(S,
T)$ is a polyhedral metric $l: E$ $\to \bold R_{>0}$ so that there
is a map, called the \it radius assignment,  \rm  $r: V$ $ \to
\bold R_{>0}$ with  $l(v v') = r(v)+ r(v')$ whenever the edge
$vv'$ has end points $v$ and $v'$.

\medskip
\noindent {\bf Theorem 3.2. (Thurston, Andreev) \it Suppose $(S,
T)$ is a triangulated closed surface.

(i)(\cite{Th}, \cite{An}) A Euclidean circle packing metric on
$(S, T)$ is determined up to isometry and scaling by the discrete
curvature $k_0$.

(ii) (\cite{Th}) A hyperbolic circle packing metric on $(S, T)$ is
determined up to  isometry by the discrete curvature $k_0$. \rm
\medskip

\noindent{\bf Colin de Verdi\`ere's Proof} (\cite{CV}). We will
consider the case (ii) for simplicity. The same argument also
works for (i) with a little care. For a circle packing metric, let
$r: V \to \bold R_{>0}$ be the radius assignment. Define $u: V \to
\bold R_{<0}$ by $u(v) =\int^{r(v)}_{\infty} \frac{1}{\sinh(s)}
ds$. The set of all circle packing metrics on $(S, T)$ is
parameterized by $\bold R_{<0}^V$ via $u$.

Recall that for a triangle of edge lengths $r_1+r_2, r_2+r_3,
r_3+r_1$ and inner angles $a_1, a_2, a_3$, Colin de Verdi\`ere
defines an action functional $F(r_1, r_2, r_3)$ by integrating the
1-form (2.6) so that
\begin{equation}\label{} \frac{\partial F}{\partial u_i} = a_i \end{equation} where $u_i
=\int^{r_i}_{\infty} \frac{1}{\sinh (s)} ds$  and proves that $F$
is strictly concave in $(u_1, u_2, u_3)$ in (2.6).  We call
$F(u_1, u_2, u_3)$ the Colin de Verdi\`ere energy of the triangle.
Define the energy $W(u)$ of $ u \in \bold R_{<0}^V$ to be the sum
of the Colin de Verdi\`ere energy of the triangles in the circle
packing metric $r$. Then by the construction, $W(u)$ is strictly
concave in $u$ due to the concavity of $F$. Furthermore, by (3.1),
$$\bigtriangledown W = 2\pi(1,....,1) -k_0$$ where $k_0$ is the
curvature of the circle packing metric. Thus, by lemma 3.1, the
map from $r$ to its curvature $k_0$ is injective. This is the
statement in (ii).

We remark that Colin de Verdi\`ere \cite{CV} also gave a very nice
proof of the existence of the circle packing metrics.

\medskip
\subsection{The work of Rivin and Leibon}

Given a polyhedral surface $(S, T)$, Rivin \cite{Ri} introduced
the  curvature $\phi_0: E \to \bold R$ sending an edge $e$ to
$2\pi-a-a'$ where $a,a'$ are the angles facing the edge.  Leibon
introduced in \cite{Le} the $\psi_0: E \to \bold R$ curvature
which sends an edge $e$ to $\frac{b+c-a+b'+c'-a'}{2}$ where $a,a'$
are the angles facing the edge $e$ and $b,b',c,c'$ are the angles
adjacent to $e$. The geometric meaning of $\phi_0$ and $\psi_0$
are related to the dihedral angles of the associated 3-dimensional
hyperbolic ideal polyhedra. See Figure 3.1.

\medskip

\epsfxsize=4.0truein\centerline{\epsfbox{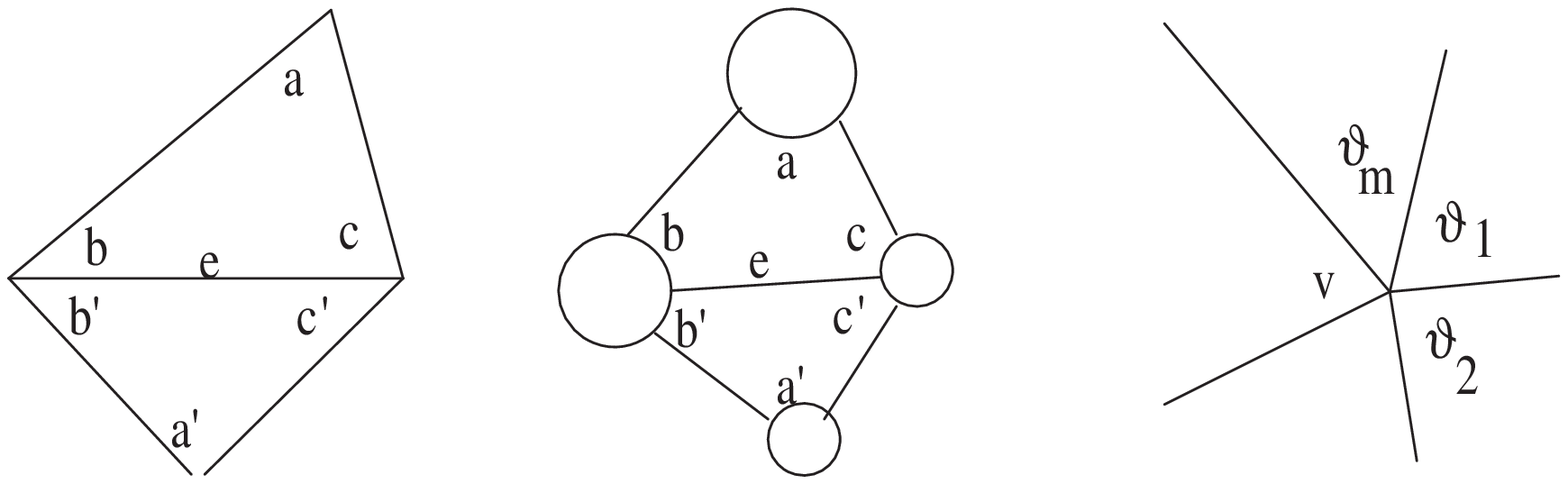}}

\centerline{Figure 3.1}

\medskip
\noindent{\bf Theorem 3.3.} (Rivin, Leibon) \it (i)(\cite{Ri}). A
Euclidean polyhedral metric on $(S, T)$ is determined up to
isometry and scaling by the $\phi_0$ curvature.

(ii) (\cite{Le}). A hyperbolic polyhedral metric on $(S, T)$ is
determined up to  isometry  by the $\psi_0$ curvature. \rm

\medskip
The proof in \cite{Ri} uses  the Lagrangian multipliers method.
The action functionals in \cite{Ri} and \cite{Le} are the
integrations of (2.8) and (2.9) and their Legendre
transformations.  Due to the convexity of the action functional,
the rigidity theorem follows essentially from lemma 3.1.

\medskip
\subsection{New curvatures for polyhedral metrics and some rigidity theorems}

Based on theorem 2.4, we introduced families of discrete
curvatures in \cite{Lu3}.  Recall that $(S, T)$ is a closed
triangulated surface so that $T$ is the triangulation, $E$ and $V$
are the sets of all edges and vertices. Let $\bold E^2$, $\bold
S^2$ and $\bold H^2$ be the Euclidean, the spherical and the
hyperbolic 2-dimensional geometries.

\medskip
\noindent {\bf Definition 3.4.} Given a $K^2$ polyhedral metric
$l$ on $(S, T)$ where $K^2$ $=\bold E^2$, or $\bold S^2$ or $\bold
H^2$, the \it $\phi_{ h  }$ curvature \rm of the polyhedral metric
$l$ is the function $\phi_{ h  }: E \to \bold R$ sending an edge
$e$ to:
\begin{equation} \label{2.1}
\phi_{ h  }(e) = \int^{\pi/2}_{a} \sin^{ h  }(t) dt +
\int^{\pi/2}_{a'} \sin^{ h  }(t) dt \end{equation} where $a, a'$
are the inner angles facing the edge $e$. See Figure 3.1.

The \it $\psi_{ h  }$ curvature \rm of the metric $l$ is the
function $\psi_{ h  }: E \to \bold R$ sending an edge $e$ to

\begin{equation}\label{2.2}
 \psi_{ h  }(e) =\int^{\frac{b+c-a}{2}}_0 \cos^{ h  }(t) dt +
\int^{\frac{b'+c'-a'}{2}}_0 \cos^{ h  }(t) dt \end{equation} where
$b,b',c,c'$ are inner angles adjacent to the edge $e$ and $a,a'$
are the angles facing the edge $e$.

The \it $ h  $-th discrete curvature $k_{ h  }$ \rm of the
polyhedral metric $l$ on $(S,T)$ is the function $k_{ h  } : V \to
\bold R$ sending a vertex $v$ to
$$k_{ h  }(v) = -\sum_{i=1}^m \int^{\theta_i}_{\pi/2}
\tan^{ h  }(t/2) dt  + (4-m)\pi/2$$ where $\theta_1$, ...,
$\theta_m$ are all inner angles at vertex $v$ and $m$ is the
degree of the vertex $v$. See Figure 3.1.

\medskip


The curvatures $\phi_0$ and $\psi_0$ were first introduced by I.
Rivin \cite{Ri} and G. Leibon \cite{Le} respectively.  The
positivity of $\psi_0$ and $\phi_0$ is shown in \cite{Ri} and
\cite{Le} to be equivalent to the Delaunay condition for
polyhedral metrics.

It is shown \cite{Lu3} that the positivity of the curvatures
$\phi_h$ and $\psi_h$ is independent of $h$. To be more precise,
due to $(x+y)( \int_0^x \cos^h(t) dt + \int^y_0 \cos^h(t) dt) \geq
0$ for $x,y \in [-\pi/2, \pi/2]$,
  we  have $\psi_h(e)
\geq 0$ (or $\phi_h(e) \geq 0$) if and only if $\psi_0(e) \geq 0$
(or $\phi_h(e) \geq 0$). Thus the geometric meaning of positive
$\psi_h$ curvature is the same Delaunay condition for polyhedral
metrics.

The curvature $\phi_{-2}(e) = \cot(a)+ \cot(a')$ has appeared  in
the finite element method approximation of the Laplace operator.
It is called the cotangent formula for discrete Laplacian.

\medskip

We prove that,
\medskip

\noindent {\bf Theorem 3.5.} (\cite{Lu3}) \it Let $  h   \in \bold
R $ and $(S, T)$ be  a closed triangulated surface.

(i)  A Euclidean circle packing metric on $(S, T)$ is determined
up to isometry and scaling by its $k_{ h }$-th discrete curvature.

(ii) A hyperbolic circle packing metric on $(S, T)$ is determined
up to  isometry by its $k_{ h  }$-th discrete curvature.

(iii) If $ h   \leq -1$,  a Euclidean polyhedral metric on $(S,
T)$ is determined up to isometry and scaling
 by its  $\phi_{ h  }$ curvature.

(iv) If $ h   \leq -1$ or $ h   \geq 0$, a spherical polyhedral
metric on $(S, T)$ is determined up to isometry by its $\phi_{ h
}$ curvature.

(v) If $ h   \leq -1$ or $ h   \geq 0$, a hyperbolic polyhedral
surface is determined up  isometry by its $\psi_{ h  }$ curvature.

 \rm
\medskip
 For any $ h   \in \bold R$,  there are
 local rigidity theorems in cases (i)-(v) (see theorem 6.2 in \cite{Lu3}).
We conjecture that theorem 3.5 holds for all $h$.
 To the best of our knowledge, theorem 3.5 for the
simplest case of the boundary of a tetrahedron is new.

The ideas of the proof are the same as the ones used in \cite{CV}
and \cite{Ri} by applying the action functionals discovered in
theorem 2.4. The extra constrains that $ h \leq -1$ or $ h  \geq
0$ in the theorem are caused by the condition on the convexity of
the domain $\Omega$ in lemma 3.1. To be more precise, these
conditions on $h$ guarantee that the corresponding domains of the
action functionals are convex.

\medskip
\subsection{Application to Teichm\"uller theory of surfaces with
boundary}

The counterpart of theorem 3.5(v) for hyperbolic metrics with
totally geodesic boundary on an ideal  triangulated compact
surface is the following. Recall that an \it ideal triangulated
compact surface \rm with boundary $(S, T)$ is obtained by removing
a small open regular neighborhood of the vertices of a
triangulation of a closed surface. The \it edges \rm  of an ideal
triangulation $T$ correspond bijectively to the edges of the
triangulation of the closed surface. Given a hyperbolic metric $l$
with geodesic boundary on an ideal  triangulated surface $(S, T)$,
there is a unique geometric ideal triangulation $T^*$ isotopic to
$T$ so that all edges are geodesics orthogonal to the boundary.
The edges in $T^*$ decompose the surface into hyperbolic
right-angled hexagons. The \it $\psi_{ h  }$ curvature \rm of the
hyperbolic metric $l$ is defined to be the map $\psi_{h}: \{$ all
edges in $T$\} $\to \bold R$ sending each edge $e$ to

\begin{equation}\label{}
\psi_{ h  }(e) = \int_0^{\frac{b+c-a}{2}} \cosh^{ h  }(t) dt
+\int^{\frac{b'+c'-a'}{2}}_0 \cosh^{ h  }(t) dt \end{equation}
where $a,a'$ are lengths of arcs in the boundary (in the ideal
triangulation $T^*$) facing the edge and $b,b',c,c',$ are the
lengths of arcs in the boundary adjacent to the edge so that
$a,b,c$ lie in a hexagon. See Figure 3.1.

\medskip
\noindent {\bf Theorem 3.6} (\cite{Lu3}). \it  A hyperbolic metric
with totally geodesic boundary on an ideal  triangulated compact
surface
 is determined up to isometry by
its $\psi_{ h  }$-curvature.  Furthermore, if $ h   \geq 0$, then
the set of all $\psi_{ h  }$ curvatures on a fixed ideal
triangulated surface is an explicit open convex polytope $P_{ h }$
in  a Euclidean space so that $P_{ h  }=P_0$. \rm

\medskip
The first part of the theorem is proved by using lemma 3.1 and the
action functional in theorem 2.4(v).

The case when $ h  <0$ has been recently established by Ren Guo
\cite{Gu1}. He proved that,

\medskip
 \noindent {\bf Theorem 3.7.} (Guo) \it Under the same
assumption as in theorem 3.6, if $ h   <0$, the set of all $\psi_{
h }$ curvatures on a fixed ideal  triangulated surface is an
explicit bounded open convex polytope $P_{ h  }$ in a Euclidean
space. Furthermore, if $ h   < \mu$, then $P_{ h  } \subset
P_{\mu}$. \rm

\medskip

 Theorem 3.6 was proved for
$h=0$ in \cite{Lu2} where the open convex polytope $P_0$ is
explicitly described. Evidently for each $ h   \in \bold R$, the
curvature $\psi_{ h  }$ can be taken as a coordinate of the
Teichm\"uller space of the surface. The interesting part of the
theorem 3.6 is that the images of the Teichm\"uller space in these
coordinates (for $ h   \geq 0$) are all the same. Whether these
coordinates are related to quantum Teichm\"uller theory is an
interesting topic. See \cite{CF}, \cite{Ka}, \cite{BL}, and others
for more information. Combining theorem 3.6 with the work of
Ushijima \cite{Us} and Kojima \cite{Ko}, one obtains for each $ h
\geq 0$ a cell decomposition of the Teichm\"uller space invariant
under the action of the mapping class group. See corollary 10.6 in
\cite{Lu3}.

\medskip
\section{The moduli spaces of polyhedral metrics}

One of the applications of the rigidity theorems is on the space
of all polyhedral surfaces. The most prominent result in the area
is the theorem of Andreev-Thurston. It states that the space of
all discrete curvatures of all circle packing metrics on a
triangulated surface is a convex polyhedron. We will present
Marden-Rodin's elegant proof \cite{MaR} of it in this section.
Another proof of it can be found in \cite{CV}. Many other results
on the space of all curvatures will also be discussed. Most of the
results obtained in \cite{Lu3} on the space of all curvatures are
modelled on the Marden-Rodin's method.

\medskip
\subsection{Thurston-Andreev's theorem and Marden-Rodin's proof}

\medskip
\noindent {\bf Theorem 4.1.} (Thurston-Andreev) (\cite{An},
\cite{Th}) \it The space of discrete curvatures $k_0$ of Euclidean
or hyperbolic circle packing metrics on a closed triangulated
surface is a convex polytope. \rm  \medskip

For simplicity, we will present  Marden-Rodin's proof of it
 for hyperbolic circle packing metrics. Let $(S , T)$ be the closed triangulated surface so that the set
of vertices is $V$.
 By definition, a circle packing metric is
given by the radius parameter $r \in \bold R^V_{>0}$, its discrete
curvature $k_0  \in \bold R^V$. By
 Thurston-Andreev's rigidity theorem 3.2, the curvature map $\Pi:
\bold R^V_{>0} \to \bold R^V$ sending a metric $r$ to its
curvature $k_0$ is a smooth embedding. The goal is to show that
the image $\Pi(\bold R^V_{>0})$ is a convex polytope. To this end,
one needs to study how hyperbolic triangles degenerate. Suppose a
hyperbolic triangle has edge lengths $r_1+r_2, r_2+r_3, r_3+r_1$
so that the angle opposite to $r_i+r_j$ is $\theta_k$,
$\{i,j,k\}=\{1,2,3\}$. We say a sequence of  hyperbolic metrics
degenerates if one of $r_i$ converges to $0$ or $\infty$.

The following simple lemma summarizes the degenerations. The best
proof of it is to draw a picture. See Figure 4.1 and \cite{MaR}
for a proof.

\medskip
\noindent {\bf Lemma 4.2.} (\cite{MaR} and \cite{Th}). \it Under
the assumption above,

\noindent(a) $\lim_{r_i \to \infty} \theta_i(r_1, r_2, r_3)=0$ so
that the convergence is uniform.

\noindent (b) Suppose $f_1, f_2, f_3$ are positive real numbers.
Then
$$\lim_{(r_i,r_j, r_k) \to (0, f_j, f_k) } \theta_i(r_1, r_2, r_3) =\pi, $$
 $$\lim_{ (r_i, r_j, r_k) \to (0, 0, f_k) } (\theta_i+\theta_j)(r_1, r_2, r_3) =\pi,$$ and
  $$\lim_{ (r_1, r_2, r_3) \to 0}(\theta_1+\theta_2+\theta_3)(r_1, r_2, r_3) =\pi.$$ \rm

\medskip

\vskip.1in

\epsfxsize=3.3truein \centerline{\epsfbox{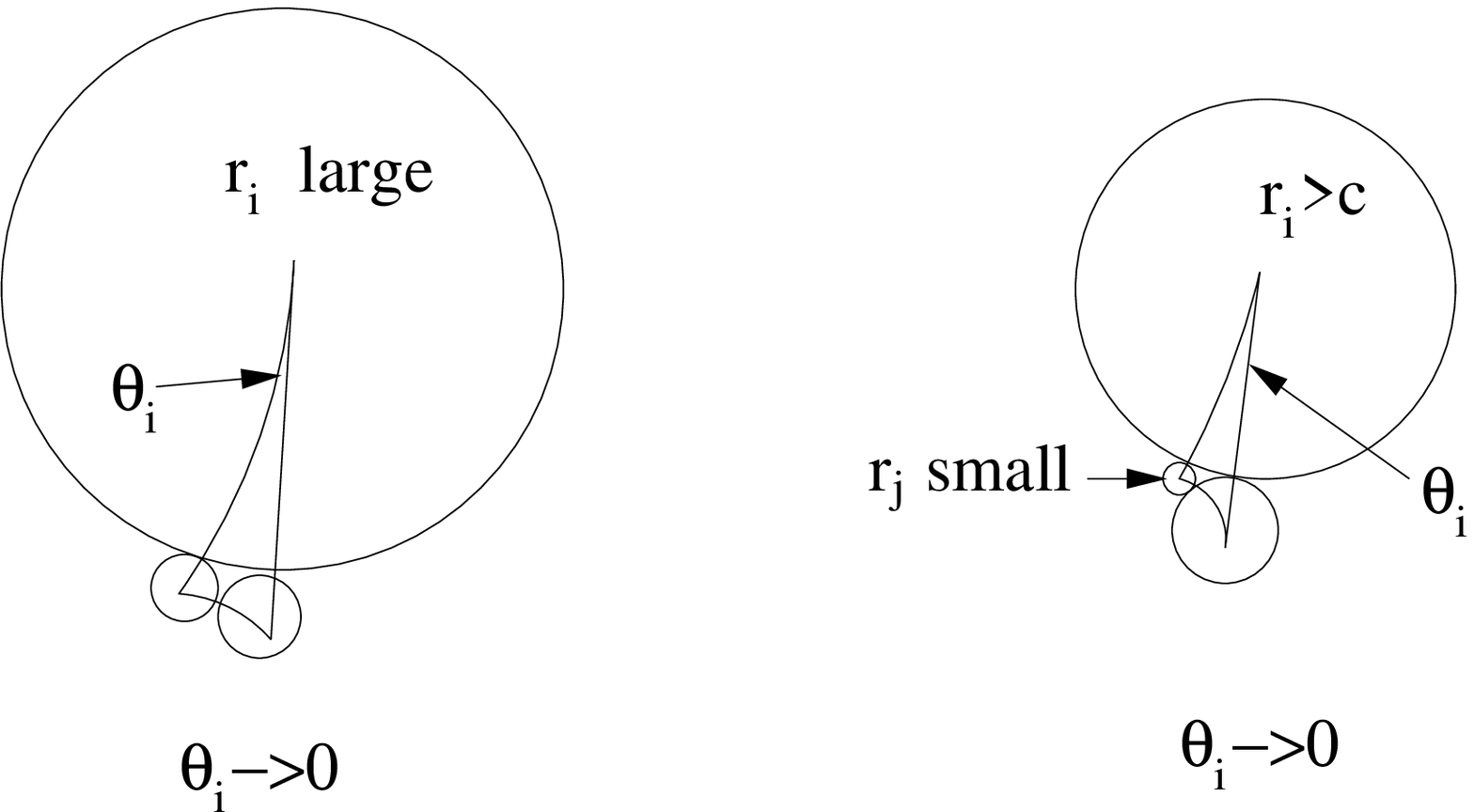}}
\centerline{Figure 4.1}



\medskip

Using the lemma, one can now determine the boundary of $\Omega
=\Pi(\bold R^V_{>0})$ as follows. By the definition of the
discrete curvature $k_0$, the image $\Omega$ lies in the open
half-spaces \begin{equation} \label{} k_0(v) < 2\pi.
\end{equation} To see that the equality holds in the limit case,
 take a sequence of circle packing
metrics $r^{(n)}$ converging to a point $a$ in the boundary of
$\bold R_{>0}^V$ in the space $[0, \infty]^V$.  If one of the
coordinate of $a$, say $a(v)$,  is infinite, then by lemma 4.2,
its curvature $k_0(v)$ converges to $2\pi$, i.e., (4.1) becomes an
equality.  Now suppose $a \in [0, \infty)^V$.  Let $ I \subset V$
be the non-empty set $\{ v \in V | a(v) =0\}$. Let $F_I$ be the
set of all triangles having at least one vertex in $I$. Take a
triangle $\sigma \in F_I$. Let the vertices of $\sigma$ be $v_1,
v_2, v_3$ so that the inner angle at $v_i$ is $\theta_i$. There
are three possibilities for $\sigma$: (i) there is only one
vertex, say $v_1$, of $\sigma$ in $I$; (ii) there are two
vertices, say $v_2, v_3$, of $\sigma$ in $I$; and (iii) all
vertices $v_1, v_2, v_3$ of $\sigma$ are in $I$. In the case (i),
by lemma 4.1, $\theta_1=\pi$. In the case (ii), $\theta_2+\theta_3
=\pi$ (note that the individual value $\theta_i$ may not be well
defined). In the case (iii), $\theta_1+\theta_2+\theta_3 =\pi.$
Consider the sum of all discrete curvatures at $I$,
$$ \sum_{ v \in I} k_0(v) = 2\pi |I| - \sum_{ \theta \in X} \theta$$
where $X$ is the set of all inner angles having vertices in $I$.
We may group these inner angles in $X$ according to the triangles
they lie and sum over  each triangle first. The result is that
$\sum_{ \theta \in X} \theta = \pi |F_I|$ by the discussion above
in the three cases (i), (ii) and (iii), i.e.,
\begin{equation} \sum_{ v \in I} k_0(v) = 2\pi |I| -\pi |F_I|
\end{equation} On the other hand, it is easy to see that, \begin{equation}\label{} \sum_{ v \in I} k_0(v) > 2\pi |I| -\pi
|F_I| \end{equation} due to the fact that the sum of inner angles
of a hyperbolic triangle is less than $\pi$. Thus the image
$\Pi(\bold R_{>0}^V)$ is the open bounded convex polytope bounded
by linear inequalities (4.1) and (4.3).  This ends the proof.

\medskip
\subsection{Some other results on the space of curvatures}

Similar results for the moduli spaces of all $\phi_{ h  }$ or
$\psi_{ h  }$ curvatures, or $k_{ h  }$ discrete curvatures on a
triangulated surfaces are obtained in \cite{Lu3} using the same
method.

\medskip
\noindent {\bf Theorem 4.3.} (\cite{Lu3}). \it Suppose $ h   \leq
-1$ and $(S, T)$ is a closed triangulated surface so that $V$ is
the set of all vertices. Then,

(i) The space of all $k_{ h  }$-discrete curvatures of Euclidean
circle packing metrics on $(S, T)$ forms a proper codimension-1
smooth submanifold in $\bold R^V$.

(ii) The space of all $k_{ h  }$-discrete curvatures of hyperbolic
circle packing metrics on $(S, T)$ is an open submanifold in
$\bold R^V$ bounded by the proper codimension-1 submanifold in
part (a). \rm

\medskip

A related theorem, proved in essentially the same way,  is the
following. Given a closed triangulated surface $(S,T)$ and
$K^2=\bf S^2, \bf H^2, \bf E^2$, let $P_{K^2}(S, T)$ be the space
of all $K^2$ polyhedral metrics on $(S, T)$ parameterized by the
edge length function. Let $\Phi_h$ (resp. $\Psi_h$): $P_{K^2}(S,
T) \to \bf R^E$ be the map sending a metric to its $\phi_h$ (resp.
$\psi_h$) curvature.

\medskip

\noindent{\bf Theorem 4.4.} (\cite{Lu3})  \it Suppose $(S, T)$ is
a closed triangulated surface so that $E$ is the set of all edges.
Let $ h   \leq -1$.

(a) The space $\Phi_{ h  }(P_{E^2}(S, T))$ is a proper smooth
codimension-1 submanifold in $\bold R^E$.

(b) The space $\Psi_{ h  }(P_{H^2}(S, T))$ is an open set bounded
by $\Phi_{ h  }(P_{E^2}(S, T))$ and a finite set of linear
inequalities.
 \rm

\medskip

Using the same argument, we can give a proof (\cite{Lu3}, \S9) of
the following results of Rivin and Leibon. Recall that a Euclidean
or hyperbolic polyhedral surface is said to have Delaunay property
if its $\psi_0$ curvature is non-negative.

\medskip
\noindent{\bf Theorem 4.5.} (Rivin, Leibon) \it Suppose $(S, T)$
is a closed triangulated surface.

(i) (\cite{Ri}) The space of all $\phi_0$-curvatures of Delaunay
Euclidean polyhedral metrics on $(S, T)$ is a convex polytope.

(ii) (\cite{Le}) The space of all $\psi_0$-curvatures of Delaunay
hyperbolic polyhedral metrics on $(S, T)$ is a convex polytope.

\rm

\medskip

\subsection{A sketch of the proof theorems 3.6 and 3.7}

Suppose $E$ is the set of all edges in an ideal triangulation of a
compact surface $S$. Let  $\Psi_h: Teich(S) \to \bold R^E$ be the
map sending a hyperbolic metric to its $\psi_h$ curvature defined
in (3.4). The first part of theorem 3.6 says $\Psi_h$ is a smooth
embedding. It is proved using the variational principle associated
to the action functional in theorem 2.4(v).  To prove the second
part of theorem 3.6 and theorem 3.7, we need to show that the
image $\Psi_h(Teich(S))$ is a convex polytope.  To this end, one
must determine the degenerations of the right-angled hexagons.
This was achieved in \cite{Lu3} and \cite{Gu1}.

The result corresponding to lemma 4.2 is the following,

\medskip
\noindent {\bf Lemma 4.6.} (\cite{Gu1}, \cite{Lu3}) \it Suppose a
hyperbolic right-angled hexagon has three non-pairwise adjacent
edge lengths $l_1, l_2, l_3$ and opposite edge lengths $\theta_1,
\theta_2, \theta_3$ so that $\theta_i = r_j + r_k$,
$\{i,j,k\}=\{1,2,3\}$. Then the following holds.

(a)  $\lim_{\theta_i \to 0} l_j(\theta_1, \theta_2, \theta_3)
=\infty$ for $j \neq i$.

(b) $\lim_{l_i \to 0} r_i(l_1, l_2, l_3) =\infty$.

(c) Suppose a sequence of hexagons satisfies that $|r_1|, |r_2|,
|r_3|$ are uniformly bounded. Then $\lim_{l_i \to \infty}
\theta_j(l) \theta_k(l) =0$ so that the convergence is uniform in
$l$. \rm

\medskip
Using this lemma, by the same analysis of boundary points as in
Marden-Rodin's proof, we establish theorem 3.6 that the image of
$\Psi_{h}(T(S))$ is an open convex polytope in $\bold R^E$
independent of the choice of $h \geq 0$. In Guo's work \cite{Gu1},
he was able to push the analysis further and obtained the result
for all $h<0$.

\medskip

\section{Several open problems}

We believe that theorem 3.5 should be true for all parameters $h$.

The space of all polyhedral metrics on $(S, T)$ carries a natural
Poisson structure. It is very interesting to know if the Poisson
structure can be expressed explicitly in $\phi_h$ and $\psi_h$
coordinates for each $h \in \bold R$. See the recent nice work of
Mondello \cite{Mo} for related works.

By taking a sequence of polyhedral metrics converging to a smooth
Riemannian metric, we would like to know if the corresponding
discrete curvatures $\phi_h$, $\psi_h$ and $k_h$ converge in the
sense of measure.

The following question concerning the metric-curvature
relationships may have an affirmative answer.

\medskip
\noindent {\bf Problem 5.1.} (\cite{Lu3}) \it Suppose $(S, T)$ is
a closed triangulated surface. Let $\Pi: P_{K^2}(S, T) \to \bold
R^V$ be the curvature map sending a metric to its discrete
curvature $k_0$. Take $p \in \bold R^V$.

(a) For $K^2=\bold E^2$ or $\bold H^2$, show that the space
$\Pi^{-1}(p)$ is either the empty set or a smooth manifold
diffeomorphic to $\bold R^{|E|-|V|}$.

(b) For $K^2=\bold S^2$, show that the space $\Pi^{-1}(p)$ is
either the empty set or a smooth manifold diffeomorphic to $\bold
R^{ |E|-|V| + \mu}$ where $\mu$ is the dimension of the group of
conformal automorphisms of a spherical polyhedral metric $l \in
\Pi^{-1}(p)$. \rm

\medskip
Given a spherical polyhedral metric $l$ on $(S, T)$, let $V'$ be
the set of all vertices so that the discrete curvatures at the
vertices are zero. The number $\mu$ above is the dimension of the
group of all conformal automorphisms of the Riemann surface $S-V'$
where the conformal structure is induced by $l$.
 In particular, if the Euler
characteristic of $ S-V'$ is negative, then $\mu=0$.

We have shown in \cite{Lu3} that the  preimage $\Pi^{-1}(p)$ is
either empty or a smooth manifold of dimension $|E|-|V|$ in the
Euclidean or hyperbolic cases.


\frenchspacing
\begin{thebibliography}{1}


\bibitem{An}
Andreev, E. M., Convex polyhedra in Lobachevsky spaces. (Russian)
Mat. Sb. (N.S.) 81 (123) 1970 445--478.


\bibitem{BS}
Bobenko, Alexander I.; Springborn, Boris A., Variational
principles for circle patterns and Koebe's theorem. Trans. Amer.
Math. Soc. 356 (2004), no. 2, 659--689.

\bibitem{BL} Bonahon, Francis; Liu, Xiaobo, Representations of the quantum Teichmüller
space and invariants of surface diffeomorphisms. Geom. Topol. 11
(2007), 889--937.


\bibitem{Br}
Br\"{a}gger, W., Kreispackungen und Triangulierungen. Enseign.
Math., 38 (1992), 200--217.



\bibitem{CF}  Chekhov, L. O. Fok, V.V.: Quantum Teichmüller spaces. (Russian)
 Teoret. Mat. Fiz. 120 (1999), no. 3, 511--528; translation in Theoret. and Math. Phys. 120 (1999), no. 3, 1245--1259

\bibitem{CV}
Colin de Verdi\`ere, Yves, Un principe variationnel pour les
empilements de cercles. Invent. Math. 104 (1991), no. 3, 655--669.

\bibitem{CKP}
Cohn, Henry; Kenyon, Richard; Propp, James, A variational
principle for domino tiling. J. Amer. Math. Soc. 14 (2001), no. 2,
297--346.

\bibitem{CL}
Chow Bennett; Luo Feng: Combinatorial ricci flows on surfaces.
Journal Differential Geometry, 63 (2003), no.1, 197--129.


\bibitem{Gu1}
Guo, Ren, On parameterizations of Teichm¡§uller spaces of surfaces
with boundary, arXiv: math.GT/0612221.


\bibitem{Ka} Kashaev, R. M. Quantization of Teichmüller spaces
 and the quantum dilogarithm. Lett. Math. Phys. 43 (1998), no. 2, 105--115.

\bibitem{Ko}
Kojima, Sadayoshi, Polyhedral decomposition of hyperbolic
3-manifolds with totally geodesic boundary. Aspects of
low-dimensional manifolds, 93--112, Adv. Stud. Pure Math., 20,
Kinokuniya, Tokyo, 1992.

\bibitem{Le}
Leibon, Gregory, Characterizing the Delaunay decompositions of
compact hyperbolic surfaces. Geom. Topol. 6 (2002), 361--391.


\bibitem{Lu1}
Luo, Feng, A characterization of spherical polyhedron surfaces, J.
Differential Geom. 74, (2006), no. 3, 407-424.

\bibitem{Lu2}
Luo, Feng, On Teichmüller spaces of surfaces with boundary, Duke
Mathematical Journal. Volume 139, Number 3 (2007), 463-482.


\bibitem{Lu3} Luo, Feng, Rigidity of polyhedral surfaces,
arXiv:math.GT/0612714.


\bibitem{MaR}
Marden, Al; Rodin, Burt, On Thurston's formulation and proof of
Andreev's theorem. Computational methods and function theory
(Valparalo, 1989), 103--115, Lecture Notes in Math., 1435,
Springer, Berlin, 1990.


\bibitem{Mo} Mondello, Gabriele;
Triangulated Riemann surfaces with boundary and the Weil-Petersson
Poisson structure, arXiv: math/0610698, to appear in J.
Differential Geom.


\bibitem{Ri}
Rivin, Igor, Euclidean structures on simplicial surfaces and
hyperbolic volume. Ann. of Math. (2) 139 (1994), no. 3, 553--580.



\bibitem{Sch}  Schlenker, Jean-Marc, Circle patterns on singular
surfaces, arXiv:math.DG/0601531

\bibitem{Sp}  A. Springborn, Boris, A variational principle for weighted
Delaunay triangulations and hyperideal polyhedra, arXiv:
math/0603097.


\bibitem{Th}
Thurston, William P. Three-dimensional geometry and topology.
 1979-1981,
http://www.msri.org/publications/books/gt3m/

\bibitem{Us}
Ushijima, Akira. A canonical cellular decomposition of the
Teichm\"uller space of compact surfaces with boundary. Comm. Math.
Phys. 201 (1999), no. 2, 305--326.




\end{thebibliography}
\end{document}